\documentclass[12pt]{article}
\usepackage{latexsym,amsfonts,amssymb}
\usepackage[dvips]{color}
\usepackage{colordvi,multicol}

\def \cal{\mathcal}

\newtheorem{thm}{Theorem}[section]
\newtheorem{cor}[thm]{Corollary}
\newtheorem{lem}[thm]{Lemma}
\newtheorem{pro}[thm]{Proposition}
\newtheorem{defi}[thm]{Definition}

\begin{document}
\title{\bf  Convergence rates in the law of large numbers under sublinear expectations}
\author{Ze-Chun Hu, Ning-Hua Liu, Ting Ma\thanks{Corresponding
author: College of Mathematics, Sichuan University, Chengdu
610065,  China\vskip 0cm E-mail address: zchu@scu.edu.cn (Z.-C. Hu), 2639428701@qq.com  (N.-H. Liu), matingting2008@scu.edu.cn (T. Ma)}\\
 {\small College of Mathematics, Sichuan University, China}}

\maketitle
\date{}

\noindent {\bf Abstract}\ \ In this note, we study convergence rates in the law of large numbers for independent and identically distributed random variables under sublinear expectations. We obtain a strong $L^p$-convergence version and a strongly quasi sure convergence version of the law of large numbers.




\section{Introduction}

Let $\{X,X_n,n\geq 1\}$ be a sequence of independent and identically distributed (i.i.d.) random variables in a probability space $(\Omega,\mathcal{F},P)$. Define $S_n=X_1+\cdots+X_n$, $n\in\mathbb{N}$. If $E|X|<\infty$, then by the law of large numbers, we know that $\frac{S_n}{n}\stackrel{a.s.}{\longrightarrow} \mu$, where $\mu=E[X]$. In fact, we also have $\frac{S_n}{n}\stackrel{L^1}{\longrightarrow} \mu$ by the martingale theory.

Hsu and Robbins (1947) introduced a new kind of convergence named ``{\it complete convergence}''. Let $\{Y,Y_n,n\geq 1\}$ be a sequence of  random variables. $\{Y_n,n\geq 1\}$ is said to completely converge to $Y$, if for any $\varepsilon>0$,
\begin{eqnarray*}
\sum_{n=1}^{\infty}P(|Y_n-Y|>\varepsilon)<\infty,
\end{eqnarray*}
which is denoted by $Y_n\stackrel{c.c.}{\longrightarrow} Y$. Obviously, $Y_n\stackrel{c.c.}{\longrightarrow} Y\Rightarrow Y_n\stackrel{P}{\longrightarrow} Y$. In fact, by  the Borel-Cantelli lemma, we know that  $Y_n\stackrel{c.c.}{\longrightarrow} Y \Rightarrow Y_n\stackrel{a.s.}{\longrightarrow} Y$.

 Hsu and Robbins (1947) proved that if $E[X^2]<\infty$ and $E[X]=\mu$, then
$\frac{S_n}{n}\stackrel{c.c.}{\longrightarrow} \mu$. Erd\"{o}s (1949) proved the converse result. Baum and Katz (1965) extended the Hsu-Robbins-Erd\"{o}s theorem. Below is
a special case of the Baum-Katz theorem.

\begin{thm} {\rm (Baum and Katz (1965))}. Let $\alpha\geq 1$. Suppose that $\{X, X_n,n\geq 1\}$ is a sequence of i.i.d.
random variables with partial sum $S_n=\sum_{i=1}^n X_i$,
$n\in\mathbb{N}$. Then the condition $E|X|^{\alpha}<\infty$ and
$EX=0$ is equivalent to
$
\sum_{n=1}^{\infty}n^{\alpha-2}P\left(|S_n|>n\varepsilon\right)<\infty,\ \
\forall\varepsilon>0.
$
\end{thm}
Lanzinger (1998), Gut and Stadtm\"{u}ller (2011),
Chen and Sung (2014) extended the results of Baum and Katz (1965).

Chow (1988) first investigated the complete moment
convergence and obtained the following result. Let $\alpha\geq 1$,
$p\le \alpha$ and $p<2$. Suppose that $\{X, X_n,n\geq 1\}$ is a
sequence of i.i.d. random variables with $E[X]=0$. If
$E[|X|^{\alpha}+|X|\log^+|X|]<\infty$, then
\begin{eqnarray*}\label{Chow-a}
\sum_{n=1}^{\infty}n^{\frac{\alpha}{p}-\frac{1}{p}-2}E\left[\left(|S_n|-\varepsilon
n^{\frac{1}{p}}\right)^{+}\right]<\infty,\ \forall
\varepsilon>0,
\end{eqnarray*}
where $x^+=\max\{0,x\}$. Chow's result has been generalized in various directions. Refer to Qiu and Chen (2014), Li and Hu (2017) and the references therein.

Li and Hu (2017) introduced a new convergence called ``{\it strong $L^p$-convergence}''. Let $\{Y,Y_n,n\geq 1\}$ be a sequence of random variables, and $p>0$. $Y_n$ is said to strongly $L^p$-converge to $Y$ if
$
\sum_{n=1}^{\infty}E[|Y_n-Y|^p]<\infty,
$
which is denoted by $Y_n\stackrel{S\mbox{-}L^p}{\longrightarrow} Y$. Obviously, $Y_n\stackrel{S\mbox{-}L^p}{\longrightarrow} Y\Rightarrow Y_n\stackrel{L^p}{\longrightarrow} Y$.
 By Markov's inequality, $Y_n\stackrel{S\mbox{-}L^{p}}{\longrightarrow} Y\Rightarrow Y_n\stackrel{c.c.}{\longrightarrow} Y$. Then for $p\geq 1$, we have the following diagram:
\begin{eqnarray*}
\begin{array}{ccccccc}
Y_n\stackrel{c.c.}{\longrightarrow} Y&\Rightarrow& Y_n\stackrel{a.s.}{\longrightarrow} Y&\Rightarrow &Y_n\stackrel{P}{\longrightarrow} Y&\Rightarrow &Y_n\stackrel{d}{\longrightarrow} Y\\
\Uparrow & & & & \Uparrow & &\\
Y_n\stackrel{S\mbox{-}L^{p}}{\longrightarrow} Y&\Rightarrow &Y_n\stackrel{L^{p}}{\longrightarrow} Y&\Rightarrow &Y_n\stackrel{L^1}{\longrightarrow} Y&&
\end{array}
\end{eqnarray*}
In particular, for $p>1$,  $Y_n\stackrel{S\mbox{-}L^{p}}{\longrightarrow} Y$ implies both $Y_n\stackrel{a.s.}{\longrightarrow} Y$  and $Y_n\stackrel{L^1}{\longrightarrow} Y$.

Recently, Hu and Sun (2018) studied convergence rates in the law of large numbers for i.i.d. random variables. They obtained a strong $L^p$-convergence version and a strong almost sure convergence version of the law of large numbers in a probability space.

The motivation of this note is to study convergence rates in the law of large numbers for i.i.d. random variables under sublinear expectations, and extend some results in a probability space to a sublinear expectation space.

Motivated by the risk measures, superhedge pricing and modeling uncertainty in finance, Peng (2006, 2007, 2008a, 2008b, 2009, 2010) initiated the notion of i.i.d. random variables under sublinear expectations, and proved the central limit theorems and the weak law of large numbers among others.

 Hu and Zhou (2015) presented some multi-dimensional laws of large numbers under sublinear expectations without the requirement of identical distribution.  Chen (2016) proved a
 strong law of large numbers (SLLNs) for i.i.d. random variables
 under capacities induced by sublinear  expectations.  Hu and Chen (2016)
 presented three laws of large numbers for independent random
 variables without the requirement of identical distribution.    Zhang (2016) showed that Kolmogorov's SLLNs holds for i.i.d. random variables under a continuous sublinear expectation if and only if the corresponding Choquet integral is finite. Chen et al. (2017)  investigated some SLLNs for sublinear expectation without independence. Hu and Yang (2017) obtained  a SLLNs for i.i.d. random variables under one-order type moment condition. Hu (2018) obtained a SLLNs for a sequence of independent random variables satisfying a controlled 1st moment condition under sublinear expectations. Chen et al. (2019) established a kind of SLLNs for capacities with a new notion of exponential independence for random variables under an upper expectation.

 We refer to Marinacci (1999), Maccheroni and Marinacci (2005), Cozman (2010), Li and Chen (2011), Chen (2012), Chen et al. (2013), Agahi et al. (2013), Zhang and Chen (2015), Hu et al. (2016), Wu and Jiang (2018) for more results on SLLNs for capacity, nonlinear expectations or sublinear expectations.  We also refer to  Hu and Zhou (2019) and Zhang (2019) for some recent results on the convergence of random variables  under sublinear expectations.

The rest of this note is organized as follows. In Section 2, we recall some basic notions and  results on sublinear expectations. In Section 3, we present our main results and give the proofs. In the final section, we mention some questions.

\section{Sublinear expectations}\setcounter{equation}{0}\label{sub-e}

  In this section, we introduce some basic definitions and notations about sublinear expectation. Refer to Peng (2010) for more details.

Let $\Omega$ be a given set and let $\mathcal{H}$ be a linear space of real functions defined on $\Omega$ such that for any constant number $c$, $c \in \mathcal{H}$; if $X \in \mathcal{H}$, then $|X| \in \mathcal{H}$; if $X_1,\cdots ,X_n \in \mathcal{H}$, then for any
$\varphi \in C_{l,lip}(\mathbb{R}^n)$, $\varphi(X_1,\cdots ,X_n) \in \mathcal{H}$, where $ C_{l,lip}(\mathbb{R}^n)$ denotes the linear space of functions $\varphi$ satisfying
 \begin{eqnarray*}
 |\varphi(x)-\varphi(y)| \leq C(1+|x|^m+|y|^m)|x-y|, \ \ \forall x, y \in\mathbb{R}^n
 \end{eqnarray*}
for some $C>0$, $m\in \mathbb{N}$ depending on $\varphi$.
For $n \in \mathbb{N}$, let $\mathcal{H}^{n}=\{X=(X_1,\cdots,X_n), X_i\in \mathcal{H}, \forall i=1,\cdots,n\}$ denote $n$-dimensional random vector space.
\begin{defi} {\rm (Definition 1.1 of Peng (2010))}  A sublinear expectation $\hat{E}$ on $\mathcal{H}$ is a functional $\hat{E}:\mathcal{H} \rightarrow \mathbb{R}$ satisfying the following properties $:\forall$ $X,Y\in \mathcal{H}$,

(i) Monotonicity: $\hat{E}[X] \geq \hat{E}[Y] ,\  \mbox{if}\ X \geq Y$;

(ii) Constant preserving: $\hat{E}[c]=c ,\  \forall \ c \in \mathbb{R}$;

(iii) Sub-additivity:  $\hat{E}[X+Y] \leq \hat{E}[X] + \hat{E}[Y]$;

(iv) Positive homogeneity:  $\hat{E}[\lambda X] = \lambda \hat{E}[X], \  \forall \lambda \geq 0$.\\
The triple $(\Omega,\mathcal{H},\hat{E})$ is called a sublinear expectation space.

\end{defi}
In the following, we assume that $\Omega$ is a complete separable metric space and let $\cal{B}(\Omega)$ denote
the Borel $\sigma$-algebra of $\Omega$. Further we assume that there exists a family $\mathcal{P}$ of probability measures on $(\Omega,\cal{B}(\Omega))$ such that
\begin{eqnarray*}
\hat{E}[X]=\sup\limits_{P\in \mathcal{P}}{E_P}[X],\ \ \forall X \in \mathcal{H}.
\end{eqnarray*}
Suppose that for any $A \in \cal{B}(\Omega)$, $I_{A}\in \mathcal{H}$. A pair of capacities associated with $\hat{E}[\cdot]$ are defined by
\begin{eqnarray*}
V(A):=\hat{E}[I_{A}],\  v(A):=-\hat{E}[-I_{A}],\ \forall A\in \cal{B}(\Omega).
\end{eqnarray*}
It is easy to check that
\begin{eqnarray*}
V(\emptyset)=0, \ V(\Omega)=1 , V(A)+v(A^{c})=1,
\end{eqnarray*}
where $A^{c}$ is the complementary set of $A$, $A\in \cal{B}(\Omega)$.
For $p\in [1,+\infty)$, the map
\begin{eqnarray*}
\left\|\cdot\right\|_{p}:X\in \mathcal{H} \mapsto (\hat{E}[|X|^{p}])^{\frac{1}{p}}
\end{eqnarray*}
forms a seminorm on $\mathcal{H}$.


\begin{defi}{\rm (Definition I.3.1 of Peng (2010))}  Let $(\Omega_i,\mathcal{H}_i,\hat{E}_i),i=1,2$ be two sublinear expectation spaces and $X_i\in \mathcal{H}_i^{n},i=1,2$. $X_1$ and $X_2$ are called identically distributed, which is denoted by $X_1\stackrel{d}{=}X_2$, if
\begin{eqnarray*}
\hat{E}_{1}[\varphi(X_1)]=\hat{E}_{2}[\varphi(X_2)] ,\forall \varphi \in C_{l,lip}(\mathbb{R}^n).
\end{eqnarray*}
\end{defi}

\begin{defi}{\rm (Definition I.3.10 of Peng (2010))}  Let $(\Omega,\mathcal{H},\hat{E})$ be a sublinear expectation space, and $X \in \mathcal{H}^m$, $Y\in \mathcal{H}^n$, $n,m\in\mathbb{N}$. $Y$ is said to be independent to  $X$ under $\hat{E}[\cdot]$, if for each test function $\varphi \in C_{l,lip}(\mathbb{R}^{n+m})$, we have
\begin{eqnarray*}
\hat{E}[\varphi(X,Y)]=\hat{E}[\hat{E}[\varphi(x,Y)]_{x=X}],
\end{eqnarray*}
whenever $\bar{\varphi}(x):=\hat{E}[|\varphi(x,Y)|]<\infty$ for all $x$ and $\hat{E}[|\bar{\varphi}(X)|]<\infty$.
\end{defi}

\begin{defi}{\rm (Proposition I.3.15 of Peng (2010))}  A sequence of random variables $\{X_n,n\geq 1\}$ on $(\Omega,\mathcal{H},\hat{E})$ is said to be independent and identically distributed, if $X_i\stackrel{d}{=}X_1$ and $X_{i+1}$ is independent to $(X_1,\cdots,X_i)$ for each $i\geq 1$.
\end{defi}

\begin{defi}{\rm (Definition II.1.4 of Peng (2010))}  A $d$-dimensional random vector $X=(X_1,\cdots,X_d)^{T}$ on a sublinear expectation space $(\Omega,\mathcal{H},\hat{E})$ is called (centralized ) G-normal distributed if
\begin{eqnarray*}
aX+b\bar{X}\stackrel{d}{=} \sqrt{a^2+b^2}X\ \ \mbox{for any}\ a,b\geq0,
\end{eqnarray*}
where $\bar{X}$ is an independent copy of $X$ ($\bar{X}\stackrel{d}{=}X$ and $\bar{X}$ is independent to $X$).
\end{defi}


\section{Convergence rates in the law of large numbers}\setcounter{equation}{0}

In this section, we will study convergence rates in the law of large numbers under sublinear expectations. Let $(\Omega,\mathcal{H},\hat{E})$ be a sublinear expectation space as introduced in Section \ref{sub-e}, and  $\{X,X_n,n\geq 1\}$ be a sequence of random variables in $\mathcal{H}$. We have the following convergences:
\begin{description}
\item[(1)] $\{X_n,n\geq 1\}$ is said to quasi surely converge to $X$, if there exists a set $N \subset \Omega$ such that $\hat{E}[I_{N}]= 0$ and $\forall \omega\in \Omega\backslash N$, $\lim_{n \to \infty}X_n(\omega)=X(\omega)$, which is denoted by $X_n\stackrel{q.s}{\longrightarrow} X$.

\item[(2)] $\{X_n,n\geq 1\}$ is said to converge to $X$ in capacity, if for all $\varepsilon > 0$, $\lim_{n \to \infty}V(\{|X_n-X|\geq \varepsilon\})=0$, which is denoted by $X_n\stackrel{V}{\longrightarrow}X$.

\item[(3)] $\{X_n,n\geq 1\}$ is said to $L^p$ converge to $X~(p >0)$, if $\lim_{n \to \infty}\hat{E}[|X_n-X|^p]=0$, which is denoted by $X_n\stackrel{L^p}{\longrightarrow}X$.

\item[(4)] $\{X_n,n\geq 1\}$ is said to completely converge to $X$, if for any $\varepsilon >0$, $\sum_{n=1}^{\infty}V(\{|X_n-X|\geq \varepsilon\})<\infty$, which is denoted by $X_n\stackrel{c.c.}{\longrightarrow}X$.

\item[(5)] $\{X_n,n\geq 1\}$ is said to $S\mbox{-}L^p$ converge to $X~(p >0)$, if $\sum_{n=1}^{\infty}\hat{E}[|X_n-X|^p]<\infty$, which is denoted by $X_n\stackrel{S\mbox{-}L^p}{\longrightarrow}X$.

\item[(6)] $\{X_n,n\geq 1\}$ is said to strongly quasi surely converge to $X$ with order $\alpha$~($\alpha >0$), if $\sum_{n=1}^{\infty}|X_n-X|^{\alpha}<\infty$ q.s.,
which is denoted by $X_n\stackrel{S_{\alpha}\mbox{-}q.s.}{\longrightarrow}X$.
\end{description}
Generally, we have
\begin{eqnarray*}
\begin{array}{ccccccccc}
&&&&X_n\stackrel{S_{p}\mbox{-}q.s.}{\longrightarrow}X & \Rightarrow & X_n\stackrel{q.s}{\longrightarrow} X&&\\
&&&&\Uparrow & &\Uparrow&&\\
X_n\stackrel{V}{\longrightarrow}X&\Leftarrow& X_n\stackrel{L^p}{\longrightarrow}X&\Leftarrow& X_n\stackrel{S\mbox{-}L^p}{\longrightarrow}X& \Rightarrow& X_n\stackrel{c.c.}{\longrightarrow}X&\Rightarrow&X_n\stackrel{V}{\longrightarrow}X\\
\end{array}
\end{eqnarray*}
If $\hat{E}$ has the monotone continuity property (\cite[Definition 2.2(vii)]{CJP11}), i.e. for any $X_n\downarrow 0$ on $\Omega$,  $\hat{E}[X_n]\downarrow 0$, then
we have (see Hu and Zhou (2019))
\begin{eqnarray*}
 X_n\stackrel{q.s.}{\longrightarrow} X \Rightarrow  X_n\stackrel{V}{\longrightarrow} X \Rightarrow  X_n\stackrel{d}{\longrightarrow} X.
\end{eqnarray*}

Let $\{X_n,n\geq 1\}$ be a sequence of i.i.d. random variables such that $\hat{E}[X_n]=-\hat{E}[-X_n]=\mu$. Denote $S_n=X_1+ X_2\cdots+X_n$, $\tilde{S}_n=\sum_{i=1}^{n}(X_i-\hat{E}[X_i])$. Then
$
\frac{S_n}{n}-\mu=\frac{\tilde{S}_n}{n}.
$

\subsection{Strong $L^p$ convergence version of the law of large numbers}

\begin{thm}\label{thm3.4}

Suppose that $\hat{E}[|X_1|^{\alpha}]<\infty$ for some $\alpha>2$. We have

(i) if $0< p\leq 2$, then
$
\frac{S_n}{n}\stackrel{S\mbox{-}L^{p}}{\nrightarrow} \mu;
$

(ii) if $2<p\leq \alpha$,  then
$
\frac{S_n}{n}\stackrel{S\mbox{-}L^p}{\longrightarrow} \mu.
$
\end{thm}

To prove Theorem \ref{thm3.4}, we need one lemma.

\begin{lem}\label{lem3.5} {(Theorem II.3.3, Lemma II.3.9 of Peng (2010), Theorem 3.2 of Hu (2011))} \\
 Let $\{X_i\}_{i=1}^{\infty}$ be a sequence of $\mathbb{R}^{d}$-valued i.i.d. random variables, satisfying $\hat{E}[X_i]=-\hat{E}[-X_i]=0$ and $\hat{E}[|X_i|^{2+\beta}]<\infty$ for some $\beta>0$. Then the sequence $\{\overline{S_n}\}$ defined by $\bar{S}_n:=(\sum_{i=1}^{n}X_i)/\sqrt{n}$ converges in law to $\xi$, i.e.,
\begin{eqnarray*}
\lim_{n\to\infty}\hat{E}[\varphi(\overline{S_n})]=\hat{E}[\varphi(\xi)],
\end{eqnarray*}
for any continuous function $\varphi\in C(\mathbb{R}^{d})$ satisfying the growth condition that $|\varphi(x)|\leq C(1+|x|^p)$ for some constants $C>0,~p>0$, where $\xi$ is  G-normal distributed with the law $N(\{0\}\times [\underline{\sigma}^2,\bar{\sigma}^2])$, $\underline{\sigma}^2=-\hat{E}[-X_i^2],\bar{\sigma}^2=\hat{E}[X_i^2]$.
\end{lem}


\noindent {\bf Proof of Theorem \ref{thm3.4}.}
Set $\tilde{X}_i=X_i-\mu$, $i\geq 1$. Then we have $\hat{E}[\tilde{X}_i]=\hat{E}[-\tilde{X}_i]=0,\forall i\geq 1$.

(i) Let $\hat{E}[\tilde{X}_i^2]=\bar{\sigma}^2$. Then $\bar{\sigma}^2>0$ by the assumption. By the  positive homogeneity of $\hat{E}$, we have
\begin{eqnarray*}
\hat{E}\left[\left|\frac{\tilde{S}_n}{n}\right|^{p}\right]
=\frac{1}{n^{p/2}}\hat{E}\left[\left|\frac{\tilde{S}_n}{\sqrt{n}}\right|^{p}\right].
\end{eqnarray*}
By Lemma \ref{lem3.5}, we have
\begin{eqnarray}\label{proof-lem3.4-a}
\hat{E}\left[\left|\frac{\tilde{S}_n}{\sqrt{n}}\right|^p\right]\to \hat{E}[|\xi|^p], \ \ \mbox{as}\ \ n\to\infty,
\end{eqnarray}
where $\xi$ is $G$-distributed with $\hat{E}[\xi^2]=\bar{\sigma}^2>0$. It follows that $V(\{\omega\in \Omega| |\xi(\omega)|>0\})>0$, which implies that $\hat{E}[|\xi|^p]>0$ for any $p>0$. Denote $c_p=\hat{E}[|\xi|^p]$.  Then by (\ref{proof-lem3.4-a}), there exists  $N\in \mathbb{N}$  such that
\begin{eqnarray*}
\hat{E}\left[\left|\frac{\tilde{S}_n}{\sqrt{n}}\right|^{p}\right]\geq\frac{c_p}{2}, \ \ \forall\ n\geq N.
\end{eqnarray*}
Therefore, for any $0< p \leq 2$, we have
\begin{eqnarray*}
\sum_{n=1}^{\infty}\hat{E}\left[\left|\frac{\tilde{S}_n}{n}\right|^p\right]\geq \sum_{n=N}^{\infty}\frac{1}{n^{p/2}}\hat{E}\left[\left|\frac{\tilde{S}_n}
{\sqrt{n}}\right|^{p}\right]\geq \frac{c_p}{2}\sum_{n=N}^{\infty}\frac{1}{n^{p/2}}=\infty.
\end{eqnarray*}

(ii) By the assumption, we know that for any $k\geq 1$, $X_{k+1}$ is independent to $(X_1,\ldots,X_k)$, which implies that $X_{k+1}$ is negatively dependent to $(X_1,\ldots,X_k)$ (see Zhang (2016, Definition 1.5) for the definition of {\it negative dependence}).  By the Marcinkiewicz-Zygmund inequality under sublinear expectations (see Zhang (2016, (2.13)),  the fact that $\hat{E}[\tilde{X}_i]=\hat{E}[-\tilde{X}_i]=0,\forall i\geq 1$, and Minkowski's inequality (see Peng (2010, Proposition I.4.2),  we have
\begin{eqnarray}\label{thm3.4-c}
\hat{E}\left[\left|\tilde{S}_n\right|^{\alpha}\right]
&\leq& \hat{E}\left[\max\limits_{k\leq n}\left|\tilde{S}_k\right|^{\alpha}\right]\nonumber\\
&\leq& C_{\alpha}\left\{\left(\sum_{k=1}^{n}\left(\left(\hat{E}\left[\tilde{X}_k\right]\right)^{+}
+\left(\hat{\mathcal{E}}\left[\tilde{X}_k\right]\right)^{-}\right)\right)^{\alpha}
+\hat{E}\left[\left(\sum_{k=1}^{n}\tilde{X_k}^2\right)^{\frac{\alpha}{2}}\right]\right\}\nonumber\\
&=&C_\alpha\hat{E}\left[\left(\tilde{X}_1^2+\cdots+\tilde{X}_n^2\right)^
{\alpha/2}\right]\nonumber\\
&=&C_\alpha\left\|\tilde{X}_1^2+\cdots+\tilde{X}_n^2\right\|_{\alpha/2}^{\alpha/2}\nonumber\\
&\leq&C_\alpha\left(\left\|\tilde{X}_1^2\right\|_{\alpha/2}+\cdots+\left\|\tilde{X}_n^2\right\|_{\alpha/2}\right)^{\alpha/2}\nonumber\\
&=&C_\alpha n^{\alpha/2}\hat{E}\left[\left|\tilde{X}_1\right|^{\alpha}\right],
\end{eqnarray}
where $x^+=\max\{0,x\}$, $x^-=\max\{0,-x\}$, and $C_\alpha$ is a positive constant depending only on $\alpha$. It follows that
\begin{eqnarray*}
\sum_{n=1}^{\infty}\hat{E}\left[\left|\frac{\tilde{S}_n}{n}\right|^{\alpha}\right]
\leq C_\alpha\sum_{n=1}^{\infty}\frac{\hat{E}[|\tilde{X}_1|^{\alpha}]}{n^{\frac{\alpha}{2}}}<\infty .
\end{eqnarray*}

For any $2<p<\alpha$, by H\"{o}lder's inequality under sublinear expectations (see Peng (2010, Proposition I.4.2)) and (\ref{thm3.4-c}),  we obtain that
\begin{eqnarray*}
\sum_{n=1}^{\infty}\hat{E}\left[\left|\frac{\tilde{S}_n}{n}\right|^p\right]
&\leq&\sum_{n=1}^{\infty}\frac{1}{n^p}\left(\hat{E}\left[\left|\tilde{S}_n\right|^{\alpha}\right]\right)^{p/\alpha}\\
&\leq&\sum_{n=1}^{\infty}\frac{1}{n^p}\left(C_\alpha n^{\alpha/2}\hat{E}\left[\left|\tilde{X}_1\right|^{\alpha}\right]\right)^{p/\alpha}\\
&=&\left(C_\alpha\hat{E}\left[\left|\tilde{X}_1\right|^{\alpha}\right]\right)^{p/\alpha}
\sum_{n=1}^{\infty}\frac{1}{n^{p/2}}<\infty.
\end{eqnarray*}
\hfill\fbox

In Chow (1988), the author also obtained the following result. Let $\{X,X_n,n\geq 1\}$ be a sequence of i.i.d. random variables with $E[X]=0$ in a probability space $(\Omega,\cal{F},P)$. Suppose that $1<\alpha<2$. If $E[|X|^{\alpha}\log^{+}|X|]<\infty$, then
\begin{eqnarray*}\label{Chow-a}
\sum_{n=1}^{\infty}n^{-2}E[|S_n|^{\alpha}]<\infty.
\end{eqnarray*}

As a consequence of Theorem \ref{thm3.4} and its proof, we obtain the following two corollaries.

\begin{cor}
Suppose that $\alpha>2$, $\hat{E}[|X_1|^{\alpha}]<\infty$, and $\hat{E}[X_1]=-\hat{E}[-X_1]=0$. Then, for any $2<p\leq\alpha$ and $\beta>(p+2)/2,$ we have
\begin{eqnarray*}
\sum_{n=1}^{\infty}n^{-\beta}\hat{E}[|S_n|^{p}]<\infty .
\end{eqnarray*}
\end{cor}

\begin{cor}
Suppose that  $X_1 \not \equiv \mu$ q.s.,  $\hat{E}[|X_1|^{\alpha}] <\infty$ for any $\alpha >0 $, and $p>0$. Then
\begin{eqnarray*}
\frac{S_n}{n}\stackrel{S-L^p}{\longrightarrow} \mu \Leftrightarrow p>2 \ .
\end{eqnarray*}
\end{cor}

\subsection{Strongly quasi sure convergence version of the law of large numbers}

\begin{pro}\label{thm3.10}
Suppose that $\hat{E}[|X_1|^{\alpha}]<\infty$ for some $\alpha>2$. Then for any
 $\beta>2$, we have
\begin{eqnarray*}
\frac{\tilde{S}_n}{n}\stackrel{S_{\beta}\mbox{-}q.s.}{\longrightarrow} 0,
\end{eqnarray*}
i.e., $\frac{S_n}{n}\stackrel{S_{\beta}\mbox{-}q.s.}{\longrightarrow} \mu$.
\end{pro}
{\bf Proof.} By Theorem \ref{thm3.4}(ii), we know that  for $2<\beta\leq \alpha$, it holds that $\frac{S_n}{n}\stackrel{S\mbox{-}L^{\beta}}{\longrightarrow} \mu,$ i.e.
$$
\sum_{k=1}^{\infty}\hat{E}\left[\left|\frac{S_n}{n}-\mu\right|^{\beta}\right]<\infty,
$$
which together with  the monotone convergence theorem (Cohen et al. (2011)) and the sublinear property implies that
\begin{eqnarray*}
\hat{E}\left[\sum_{k=1}^{\infty}\left|\frac{S_n}{n}-\mu\right|^{\beta}\right]
&=&\lim_{m\to\infty}\hat{E}\left[\sum_{k=1}^m\left|\frac{S_n}{n}-\mu\right|^{\beta}\right]\\ &\leq&\lim_{m\to\infty}\sum_{k=1}^m\hat{E}\left[\left|\frac{S_n}{n}-\mu\right|^{\beta}\right]\\
&=&\sum_{k=1}^{\infty}\hat{E}\left[\left|\frac{S_n}{n}-\mu\right|^{\beta}\right]<\infty.
\end{eqnarray*}
It follows that
\begin{eqnarray}\label{thm3.10-a}
\sum_{k=1}^{\infty}\left|\frac{S_n}{n}-\mu\right|^{\beta}<\infty,\ q.s.\ \ \mbox{for any}\ \  2<\beta\leq \alpha.
\end{eqnarray}

By the strong law of large numbers (see Theorem 1 of Chen (2016)), there exists a set $N\subset \Omega$ such that $V(N)=0$ and for any $\omega\in \Omega\backslash N$, there exists $M(\omega)\in \mathbb{N}$ such that for any $n\geq M(\omega)$,
\begin{eqnarray*}\label{thm3.10-b}
\left|\frac{S_n}{n}-\mu\right|<1.
\end{eqnarray*}
It follows that  for $\beta>\alpha$ and $\omega\in \Omega\backslash N$,
$$
\sum_{n=M(\omega)}^{\infty}\left|\frac{S_n}{n}-\mu\right|^{\beta}\leq
\sum_{n=M(\omega)}^{\infty}\left|\frac{S_n}{n}-\mu\right|^{\alpha},
$$
which together with (\ref{thm3.10-a}) implies that for any $\beta>\alpha$,
$$
\sum_{k=1}^{\infty}\left|\frac{S_n}{n}-\mu\right|^{\beta}<\infty,\ q.s.
$$
\hfill\fbox

\section{Some questions}\setcounter{equation}{0}

In this section, we present some questions for further research.

\noindent {\bf Question 1.} Can we extend the Hsu-Robbins theorem from a probability space to a sublinear expectation space?

\noindent {\bf Question 2.} If the answer to Question 1 is affirmative, can we prove the converse result?


In fact, we can ask more questions. As to the results on the convergence rates in the law of large numbers in a probability space, we can consider the corresponding questions in a sublinear expectation space.

%
%
%

\bigskip

{ \noindent {\bf\large Acknowledgments} \quad   This work was
supported by National Natural Science Foundation of China (Grant
No. 11771309, 11871184).}

\end{document}